\documentclass[12pt,a4paper]{article}

\usepackage{setspace}

\newtheorem{defn}{Definition}[section]
\newtheorem{thm}[defn]{Theorem}

\newtheorem{lem}[defn]{Lemma}
\newtheorem{cor}[defn]{Corollary}

\newtheorem{rem}[defn]{Remark}

\newcommand\pf{\noindent{\bf Proof. }}
\newcommand\qed{\hfill\framebox(5,5){~}}

\newcommand\eps{\varepsilon}
\newcommand\vphi{\varphi}

\newcommand{\gen}[1]{\langle #1 \rangle}

\newcommand{\dents}{{\cal{D}}}
\newcommand{\Edents}{{\cal{E}}}
\newcommand{\Fdents}{{\cal{F}}}
\newcommand\mapto{\mapsto}
\newcommand{\map}{\rightarrow}
\renewcommand{\bar}{\overline}
\renewcommand\subset\subseteq
\renewcommand{\hat}{\widehat}
\newcommand{\stab}{{\textnormal{Stab}}}
\newcommand{\bistar}{\star}
\newcommand{\Aut}{\rm{Aut}}
\newcommand{\Inn}{{\rm{Inn}}}
\newcommand{\Out}{{\rm{Out}}}

\title{A class of $2$-groups admitting an action of the symmetric group of degree $3$}
\author{K. Roberts and S. Shpectorov}

\begin{document}
\maketitle

\begin{abstract}
A biextraspecial group of rank $m$ is an extension of a special
$2$-group $Q$ of the form $2^{2 + 2m}$ by $L\cong L_2(2)$, such that
the $3$-element from $L$ acts on $Q$ fixed-point-freely. Subgroups 
of this type appear in at least the sporadic simple groups $J_2$, 
$J_3$, $McL$, $Suz$, and $Co_1$. In this paper we completely classify 
biextraspecial groups, namely, we show that the rank $m$ must be 
even and for each such $m$ there exist exactly two biextraspecial 
groups $B^\eps(m)$ up to isomorphism where $\eps\in\{+,-\}$. We 
also prove that $\Out(B^\eps(m))$ is an extension of the 
$m$-dimensional orthogonal $GF(2)$-space of type $\eps$ by the 
corresponding orthogonal group. The extension is non-split except 
in a few small cases. 
\end{abstract}

\section{Introduction}

The sporadic simple groups and their properties remain a focal point 
in finite group theory. When one looks through the lists of subgroups, 
particularly 2-local subgroups, of sporadic groups, certain shapes of 
subgroups make frequent appearance. It is an interesting problem to try 
to classify the whole series of such 2-local groups and to understand 
why only finitely many of them lead to simple groups. Furthermore, 
studying properties of such groups, particularly, the automorphism 
group can help in the study of related amalgams and geometries.

This approach was taken by Ivanov and Shpectorov in \cite{is}, where 
they studied what they called tri-extraspecial groups. A tri-extraspecial 
group is defined to be a split extension of a special $2$-group 
$Q\cong 2^{3+3m}$ by $L\cong L_3(2)$ such that $Z = Z(Q)\cong 2^3$ is 
the dual natural $L$-module and $Q/Z$ has only natural $L$-modules as 
composition factors. Such groups appear prominently in at least two 
sporadic groups, $J_4$ and $Fi_{24}$. The main result of \cite{is} is 
that for each even integer $m$ and type $\eps\in\{+,-\}$ there exists a 
unique tri-extraspecial group $T^\eps(m)$. Furthermore, the outer 
automorphism group of $T^\eps(m)$ is determined. 

If one wanted to make a straightforward generalisation of 
tri-extraspecial groups by replacing the $3$ by an $n$ in the definition, 
then one would find that such groups do not exist for $n\ge 4$. Therefore 
the cases for $n=1,2,3$ are the exceptional cases for which the definition 
is meaningful. If $n=1$, we get the usual extraspecial $2$-groups and 
their classification is well documented (for example, see \cite{go} and 
\cite{dh}). It states, just like for the tri-extraspecial groups, that 
for each even integer $m$ there are exactly two extraspecial $2$-groups 
$E^\eps(m)$ up to isomorphism. In fact, this similarity was the reason 
for the name `tri-extraspecial'. The outer automorphism groups of 
extraspecial $2$-groups can be found in \cite{g}. Therefore, to complete 
the list of possible cases we need only to study the case of $n=2$. We 
call these groups biextraspecial groups.

Thus, a biextraspecial group $G$ is defined as an extension of a $2$-group
$Q\cong 2^{2+2m}$ by $L\cong L_2(2)$ such that the centre $Z$ of $Q$ is 
the natural $L$-module and that $Q/Z$ only has natural $L$-modules as 
composition factors. We drop the condition that $G$ splits over $Q$, as 
this is automatic, and we replace the dual of the natural module with the 
natural module, since they are isomorphic as $L$-modules. 

Biextraspecial groups appear prominently among the subgroups of sporadic 
groups. We found such subgroups, for various values of $m$, in at least 
five sporadic groups: $J_2$, $J_3$, $McL$, $Suz$, and $Co_1$.

Our main results are as follows.

\begin{thm}
The rank $m$ of a biextraspecial group is even. For each even $m$, there 
are exactly two biextraspecial groups $B^\eps(m)$ up to isomorphism, 
where $\eps\in\{+,-\}$. Furthermore, $\Out(B^\eps(m))$ is an extension 
of the orthogonal $GF(2)$ space of type $\eps$ by the corresponding 
orthogonal group $O^\eps_m(2)$. The extension is non-split except when 
$m=2$ or $4$. 
\end{thm}

In Section \ref{sec:biextraspecial_groups} we state our main definition 
more precisely and show that $G$ exhibits similar structural properties 
as tri-extraspecial groups. In Section \ref{sec:dents} and 
\ref{sec:dent_space} we focus on one of the common notions between
biextraspecial and tri-extraspecial groups, namely, a related orthogonal 
space called the dent space. Following the ideology of \cite{is}, we 
define a dent to be a normal subgroup $D$ of $G$ that lies between $Q$ 
and $Z$ such that $D/Z$ is the natural $L$-module. Each pair of dents $D_1$ 
and $D_2$ has a unique diagonal dent which is denoted by $D_1+D_2$ and this 
leads to the definition of the dent spaceis an $m$-dimensional vector 
space over $GF(2)$ with respect to the addition as above. Furthermore, 
we introduce on $\dents$ a nondegenerate alternating form and a compatible 
quadratic form in terms of certain properties of dents. One consequence of 
this is that $\dents$ has even dimension, that is, $m$ is even. The type 
$\eps\in\{+,-\}$ of the orthogonal space $\dents$ becomes the type of the 
biextraspecial group $G$.  Note that here we have a difference from \cite{is},
where the quadratic form is not defined globally.

After this point we part ways with \cite{is} completely. In Section 
\ref{sec:decom} we explore the relationship between the decomposition of 
$\dents$ and $G$. We show that there is bijective correspondence between 
the decomposition of $\dents$ into the orthogonal sum of nondegenerate 
subspaces and the decomposition of $G$ into biextraspecial groups of smaller 
rank. In particular, a nondegenerate summand of dimension $k$ of $\dents$ 
corresponds to a biextraspecial subgroup of rank $k$. We then introduce 
the reverse operation called composition and denote it by $\bistar$. The 
purpose of this operation is to build biextraspecial groups from smaller 
biextraspecial pieces with the expectation that all such groups can be 
constructed from a few small examples. This is not a new idea and is a 
recurring theme in the study of extraspecial $p$-groups, as found in 
\cite{dh}, and many variations of extraspecial $p$-groups such as 
semi-extraspecial groups, as found in \cite{b}. From the standard theory 
of orthogonal spaces we know that $\dents$ decomposes into the orthogonal 
sum of $2$-dimensional orthogonal spaces and this corresponds to the 
decomposition of $G$ into biextraspecial groups of rank $2$. Thus to 
classify the biextraspecial groups we study the rank $2$ groups which we 
return to in Section \ref{sec:class}.

Before we classify the biextraspecial groups we take a detour. In Section 
\ref{sec:extrasp} we show that $\dents$ can be realised as the factor group 
of the extraspecial $2$-group $R = C_Q(t)$ where $t$ is an involution of 
$L$. Like all extraspecial $2$-groups, the factor group $\bar R = R/Z(R)$ 
has the structure of an orthogonal space. We show that $\dents$ and $\bar R$ 
are isometric orthogonal spaces. By invoking the classification of 
extraspecial $2$-groups one can show that the decompositions of 
extraspecial $2$-group $R$ and the biextraspecial $G$ match via the 
aforementioned isometry. We briefly discuss this interesting correspondence 
at the end of the section but omit the details. It should be noted that $R$ 
is a specific extraspecial $2$-group and it is not at all obvious that there 
is a one-to-one correspondence between extraspecial $2$-groups and 
biextraspecial groups. 

In Section \ref{sec:class} we classify all biextraspecial groups $G$. That is, 
for each even integer $m$ and $\eps\in\{-, +\}$ we show that there exists a 
unique biextraspecial group of rank $m$ and type $\eps$. We achieve this by 
showing that the statement is true when $m=2$ and then exploiting our 
operations of composition and decomposition.

Note that in general the normal 2-subgroup $Q$ of $G$ satisfies the properties 
of the semi-extraspecial group as defined by Beisiegel in \cite{b}. That paper, 
along with \cite{pr} by Parker and Rowley, seemed to suggest that $Q$ would 
decompose in two possible ways and this would provide us with the two types of 
biextraspecial groups. Let $S$ be a Sylow $2$-subgroup of $\rm{SL_3(4)}$ and 
$T$ be a Sylow $2$-subgroup of $\rm{SU_3(4)}$. The two likely decompositions 
of $Q$ were the central product of $n$ copies of $S$, and central product of 
$n-1$ copies of $S$ followed by a copy of $T$. Interestingly, this assertion  
turned out to be false. In both groups $B^+(m)$ and $B^-(m)$ the group $Q$ is 
in fact the same, the product of $m$ copies of $S$. 

In the final section of this paper we compute the outer automorphism group of 
$G$, $\Out(G)$. The group $\Out(G)$ can be realised as the group of automorphisms 
of $G$ that centralise $L$ and it acts on the dent space. We show that this 
action has kernel $\dents$ and its image is the full orthogonal group 
$O^\eps_m(2)$. A result from \cite{g} tells that $\Out(G)$ splits over 
$\dents$ if and only if $m=2,4$.

\section{Biextraspecial groups}
\label{sec:biextraspecial_groups}

\begin{defn}  A group $G$ is said to be a \emph{biextraspecial group
of rank $m$} if the following properties are satisfied:
\begin{itemize}
\item[{\rm{(B1)}}] $G$ is an extension of a $2$-group $Q$ of order
$2^{2+2m}$ by $L\cong L_2(2)$,
\item[{\rm{(B2)}}] The centre of $Q$, denoted by $Z$, is the natural
$L$-module (and so $Z\cong 2^2$).
\item[{\rm{(B3)}}] $\bar Q=Q/Z\cong 2^{2m}$ is elementary abelian;
furthermore, as an $L$-module, $Q/Z$ only has natural modules as
composition factors.
\end{itemize}
\end{defn}

Note that $Q$ is special and also that $Z$ and $Q$ are
characteristic in $G$. Furthermore, the second part of (B3) is
equivalent to the condition that every element of order $3$ from $G$
acts on $Q$ fixed-point-freely. As a consequence of this, we can
deduce that $G$ splits over $Q$.

\begin{lem} \label{uniquecompl} The group $G$ has a unique class of complements to $Q$.
\end{lem}

\pf Let $H = QS$ where $S$ is a Sylow $3$-subgroup of $G$. The
subgroup $H$ has index $2$ in $G$ and thus is a normal subgroup of
$G$. By the Frattini argument, $G = HN_G(S) = QN_G(S)$. Observe that
$N_Q(S) = C_Q(S) = 1$, since $S$ acts on $Q$ fixed-point-freely.
Therefore $Q\cap N_G(S) = 1$ and so $N_G(S)\cong L\cong L_2(2)$ is a
complement to $Q$ in $G$. Since every complement must have a normal
subgroup of order $3$, all complements are conjugate to
$N_G(S)$.\qed

\medskip
We now assume that $L$ is a complement to $Q$ in $G$ and write
$G = Q:L$. We note for the later part of the paper, where we compute
$\Out(G)$, that $L$ is self-normalising and that the centre of $G$
is trivial.

We now turn our attention to the structure of the $L$-module $\bar
Q$.

\begin{lem} \label{ext-VbyV} Let $V$ be an $L$-module containing a
submodule $W$ such that both $W$ and $V/W$ are isomorphic to the
natural $L$-module. Then $V$ is completely reducible, that is, it is
isomorphic to the direct sum of two natural modules.
\end{lem}

\pf Let $s,t\in L$ be of order 3 and 2, respectively. Since $\dim
C_V(t)\ge 2$, we have that $C_V(t)$ is not contained in $W$. Let
$u\in C_V(t)\setminus W$ and define $U = \gen{u, u^s, u^{s^2}}$.
Note that $s$ acts on $V$ fixed-point-freely, which means that
$u+u^s+u^{s^2}$ must be zero. It follows that $\dim U=2$. Also,
since $t$ fixes $u$, we have that $u^L=\{u,u^s,u^{s^2}\}$, which
yields that $U$ is $L$-invariant. Since both $W$ and $U$ are
irreducible, they intersect trivially. Thus, $V$ is the direct sum
of $W$ and $U$.\qed

\medskip
By induction, every $L$-module whose composition factors are all
natural modules, is completely reducible. In particular, we now have
the exact structure of $\bar Q$.

\begin{cor} \label{allisDent}
The $L$-module $\bar Q$ is the direct sum of $m$ copies of the
natural $L$-module.\qed
\end{cor}

This in turn implies the following.

\begin{cor} \label{num_of_irrs}
Every irreducible $L$-submodule of $\bar Q$ is isomorphic to the
natural $L$-module and there are exactly $2^m-1$ irreducible
submodules of $\bar Q$.\qed
\end{cor}

%--------- End of Section ---------------%

\section{Dents}
\label{sec:dents}
\begin{defn} A \emph{dent} in a biextraspecial group $G=Q:L$ is
a normal subgroup $D$ of $G$ such that $Z < D < Q$ and $D/Z$ is an
irreducible $L$-submodule of $Q/Z$.
\end{defn}

Note that any two distinct dents intersect in $Z$. Also, it follows
from Corollary \ref{num_of_irrs} that there are exactly $2^m-1$
dents.

\begin{lem} \label{abelian}
Each dent $D$ is abelian. Moreover, $D$ is isomorphic either to
$2^4$ or $4^2$.
\end{lem}

\pf Note that $\bar D \cong 2^2$. Therefore $\bar D = \gen{\bar x,
\bar y}$ for some $x,y\in D$ and, consequently, $D = \gen{x,y,Z}$.
This means that $[D,D] = \gen{[x,y]}$ is cyclic and it is contained
in $Z$. As $[D,D]$ is normal in $G$, we conclude that $[D,D]=1$, as
$Z$ is irreducible as an $L$-module. In particular, $D$ is abelian.

Let $\phi: \bar D \map Z$ be a map given by $\bar d\mapto d^2$. One
can check that $\phi$ is linear and moreover, it commutes with the
action of $L$. That is, $\phi$ is a module homomorphism between two irreducible
$L$-modules. This gives us that $\phi$ is either the trivial map or
an isomorphism. If $\phi$ is the zero map, then each element of $D$
has order $2$ and $D$ is elementary abelian $2^4$. Suppose next that
$\phi$ is invertible. Then every element of $D$ outside $Z$ has
order 4. Let again $\bar D = \gen{\bar x, \bar y}$. Then $x^2$ and
$y^2$ generate $Z$ and hence $D = \gen{x,y,Z}=\gen{x,y}\cong
4^2$.\qed

\medskip
If a dent $D$ is isomorphic to $2^4$, then we say that $D$ is
\emph{singular} and otherwise it is \emph{non-singular}. We next
determine the action of $L$ on each type of dent and show that the
action is unique for each type. Let again $L=\gen{s,t}$, where $s$
is an element of order 3 and $t$ is an involution.

For the singular dent the action is clear. Indeed, by Lemma
\ref{ext-VbyV}, the singular dent $D$ is $Z\times U$, where $U\cong Z$
 is $L$-invariant. So we can choose $a\in Z$ and $x\in U$ fixed
by $t$ and also set $b=a^s$ and $y=x^s$. Then $a$, $b$, $x$, and $y$
generate $D$ and we can compute the entire action:
\begin{equation} \label{singular_action}
\begin{array}{ccccccc}
a^t = a, & & b^t = ab, & & a^s = b, & & b^s = ab,\\
x^t = x, & & y^t = xy, & & x^s = y, & & y^s = xy.
\end{array}
\end{equation}
where $c = ab$. Such a generating set and action given in $(\ref{singular_action})$ is called the \emph{standard basis} with respect to the \emph{standard action} $L$. Note that this depends on the choice of $L$. The group $Z$ has exactly two complements in $D$ that satisfy (\ref{singular_action}); they are $\gen{x,y}$ and $\gen{ax, by}$. Indeed, $t$ fixes the unique non-identity coset $\bar x$ of $\bar D$ and $\bar x = \{ x, ax, bx, cx\}$. The elements $x$ and $ax$ are fixed by $t$ and $bx$ and $cx$ are interchanged by $t$. So indeed, there is a unique non-trivial automorphism of $D$ that commutes with the action of $L$, namely, the map that sends $x$ to $ax$ and $y$ to $by$.

Suppose now that $D$ is a non-singular dent. We claim that the generators $x$ and $y$ of $D$ can be chosen so that
\begin{equation} \label{nonsingular_action}
\begin{array}{ccccccc}
x^t = x, & & y^t = x^{-1} y^{-1}, & &  x^s = y, & & y^s = x^{-1}y^{-1}
\end{array}
\end{equation}
Indeed, let $d\in D$ so that $\bar d$ is not fixed by $t$. Set $x = dd^t$. Then $\bar x = \bar d \bar d^t \ne 1$ and hence $x$ has order $4$. Clearly, $x$ is fixed by $t$. Set $y = x^s$. Then $y$ also has order $4$. If we set $a = x^2$ and $b = y^2$, then $t$ fixes $a$ but does not fix $b$. Hence $\gen{x}\cap\gen{y} = 1$ and so $D = \gen{x,y}$. Note that $s$ fixes only the identity element of $D$ and so $xx^sx^{s^2} =1$. In particular, $x^{s^{-1}} = x^{-1}(x^{-1})^s = x^{-1} y^{-1}$ and $y^t = x^{st} = x^{ts^{-1}} = x^{-1}y^{-1}$ as desired. Again, such a generating set and action given in $(\ref{nonsingular_action})$ is called the \emph{standard basis} with respect to the \emph{standard action} $L$. Note that there are only two choices of generating set for $D$ that satisfy (\ref{nonsingular_action}); they are $\gen{x,y}$ and $\gen{x^{-1}, y^{-1}}$. Indeed, $t$ fixes the unique non-identity coset $\bar x$ of $\bar D$ and $\bar x = \{ x, x^{-1}, xy^2, x^{-1} y^2 \}$. The elements $x$ and $x^{-1}$ are fixed by $t$ and  $xy^2$ and $x^{-1} y^2$ are interchanged by $t$. In particular, there is a unique non-trivial automorphism of $D$ that commutes with the action of $L$, namely, the map that sends $x$ to $x^{-1}$ and $y$ to $y^{-1}$. We state the following result that will be used in Section \ref{sec:out}

\begin{lem} \label{uniqueauto} Every dent $D$ has a unique non-trivial automorphism that commutes with the action of $L$.\qed
\end{lem}

Note that in the non-singular case, if we write $a = x^2$ and $b = y^2$, then the action of $L$ on both types of dents is essentially given by the same formulae. In both cases, we have that $t = t_a$ is the unique involution of $L$ that fixes both $x$ and $a$. The group $L$ has two other involutions, one of which fixes $y$ and $b$ and the other fixes $w = x^{-1} y^{-1}$ and $c = ab$, which we denote by $t_b$ and $t_c$, respectively.

We now determine the commutator table of any two dents. For the rest of the section we assume that $D_i = \gen{x_i, y_i, Z}$ and $D = \gen{x,y,Z}$ are chosen so that $x_i$, $y_i$, $x$ and $y$ are part of the standard basis with respect to the action of $L$.

\begin{lem} Two dents $D_1$ and $D_2$ either commute or have a uniquely determined commutator table.
\end{lem}

\pf The commutator subgroup $[D_1, D_2]\le Z$ is generated by $[x_1,x_2]$, $[x_1, y_2]$, $[y_1, x_2]$ and $[y_1, y_2]$. The involution $t$ fixes $[x_1, x_2]$ since it fixes $x_1$ and $x_2$. Thus $[x_1, x_2]\in\{1,a\}$. Similarly, $[y_1,y_2]\in\{1,b\}$ and $[w_1, w_2]\in\{1,c\}$. Suppose that $[x_1, x_2] = 1$. Then $1 = [x_1, x_2]^s = [y_1, y_2]$ and $1 = [x_1, x_2]^{s^{-1}} = [w_1, w_2]$. Using the bilinearity of the commutator bracket, we deduce that $1 = [w_1, w_2] = [x_1 y_1, x_2 x_2] = [x_1, x_2][x_1, y_2][y_1, x_2][y_1, y_2] = [x_1, y_2][y_1, x_2]$. In particular, $[x_1, y_2] = [y_1, x_2]$. Therefore $[D_1, D_2]$ is a proper $L$-submodule of $Z$ and hence the trivial space. That is, $D_1$ and $D_2$ commute. Suppose now that $D_1$ and $D_2$ do not commute. Then $[x_1, x_2] = a$ and $[y_1, y_2] = b$. Again, we deduce that $[x_1, y_2] = [y_1, x_2]$. Note that $[x_1, y_2]^{t_c} = [y_1, x_2] = [x_1, y_2]$ and so $[x_1, y_2]\in\{1,c\}$. Similarly, $[x_1, w_2]\in\{1,b\}$. Using bilinearity we have that $a = [x_1, x_2] = [x_1, y_2][x_1, w_2]$. This implies that $[x_1, y_2] = c$. Thus, we have shown that the commutator table of $D_1$ and $D_2$ is uniquely determined.\qed

\medskip
Therefore, if $D_1$ and $D_2$ do not commute, then the table of commutators in $[D_1, D_2]$ determined by

\doublespacing
\[
\begin{tabular}{| c | c  c |}
\hline
 $[ , ]$& $x_2$ & $y_2$ \\[0.1cm]
\hline
$x_1$ & a & c \\
$y_1$ & c & b \\
\hline
\end{tabular}
\]
\singlespace

The next result follows implicitly from the proof of the commutator table above.

\begin{cor} \label{all-or-nothing}
Let $D_1$ and $D_2$ be dents that do not commute. Let $d_1\in
D_1\setminus Z$ and $d_2\in D_2\setminus Z$ such that $\stab_L(\bar d_i)
= \stab_L(z_i)$ for some $z_i\in Z$ for $i=1,2$.  Then
\begin{itemize}
\item[{\rm{(i)}}]  if $z_1 = z_2$, then $[d_1, d_2] = z_1$, and
\item[{\rm{(ii)}}] if $z_1 \ne z_2$, then $[d_1, d_2] = z_1 z_2$.
\end{itemize}
\noindent In particular, $[d_1, d_2]\neq 1$.\qed
\end{cor}

%--------------- End of Section -----------------%

\section{Dent space}
\label{sec:dent_space}
Let $\vphi$ be the unique $L$-isomorphism between $\bar D_1$ and $\bar D_2$. Define $\bar D_1 + \bar D_2$ to be the diagonal of $\bar D_1$ and $\bar D_2$, that is,
\[ \bar D_1 + \bar D_2 = \{ \bar d + \bar d\vphi \mid \bar d\in \bar D_1\}. \]
\noindent When $D_1$ and $D_2$ are distinct, we define $D_1 + D_2$ as the full preimage of $\bar D_1 + \bar D_2$ in $Q$. The group $D_1 + D_2$ is necessarily a dent because $\bar D_1 + \bar D_2$ is an irreducible $L$-submodule of $\bar Q$ and is called the diagonal dent of $D_1$ and $D_2$. If $D_1 = D_2$, then we define $D_1 + D_2$ to be the formal zero and denote it by $0_\dents$. We define $\dents$ to be the union of the set of all dents and $0_\dents$. One can show that $(\dents, +)$ forms an $m$-dimensional vector space over $GF(2)$ with additive identity $0_\dents$. We show this in Section \ref{sec:extrasp}.

We now describe the relationship between dents and the diagonal dents. This enables us to define forms on the dent space.

\begin{lem} \label{all-or-one2} Let $D_1$ and $D_2$ be distinct dents and $D_3$ be the corresponding diagonal dent. Then for a dent $D$ either
\begin{itemize}
\item[{\rm{(i)}}] all three dents commute with $D$, or
\item[{\rm{(ii)}}] exactly one dent commutes with $D$.
\end{itemize}
\end{lem}

\pf Suppose that $D_1$ and $D_2$ commute with $D$. It is clear that $D_3\le D_1 D_2$ commutes with $D$. Suppose that $D_1$ and $D_2$ do not commute with $D$. Then $[x_i, x]=a$, for $i=1,2$, and $[x_1x_2, x]=[x_1, x][x_2, x]=a^2 = 1$ where $1\ne x_1x_2 \in D_3$. Corollary \ref{all-or-nothing} implies that $[D_3, D]=1$ as required. \qed

\medskip

Lemma \ref{all-or-one2} provides us with a symplectic bilinear form on the dent space $\dents$. Namely, $\beta: \dents\times\dents\map GF(2)$ which is given by
\[ \beta(D_1, D_2) = \left\{ \begin{array}{ll} 0 & \mbox{if } [D_1, D_2]=1 \\ 1 & \mbox{otherwise} \end{array} \right. \]
Note that $\beta$ is nondegenerate since no dent lies in the centre of $Q$. As a result, the dimension of $\dents$ is even, that is, $m$ is even.

\begin{lem} \label{one-or-all} Let $D_1$ and $D_2$ be dents and $D_3$ be the corresponding diagonal dent, then
\begin{itemize}
\item[{\rm{(i)}}] if $[D_1, D_2] = 1$, then either all three dents are singular or exactly one dent is singular; and
\item[{\rm{(ii)}}] if $[D_1, D_2] \ne 1$, then either all three dents are non-singular  or exactly one dent is non-singular.
\end{itemize}
\end{lem}

\pf We start by noting that each element of a non-singular dent outside of $Z$ has order $4$. For (i) suppose that $D_1$ and $D_2$ commute. If $D_1$ and $D_2$ are singular, then $(x_1x_2)^2 = x_1 x_2 x_1 x_2 = [x_1, x_2] = 1$. The element $x_1 x_2\in D_3\setminus Z$ and so $D_3$ is singular. If $D_1$ and $D_2$ are non-singular, then $(x_1 x_2)^2 = x_1^2 x_1^2 = a^2 = 1$. Thus, $D_3$ is singular. For (ii) suppose now that $D_1$ and $D_2$ do no commute. If $D_1$ and $D_2$ are non-singular, then multiplying $(x_1x_2)^2$ by $x_1^2 x_2^2 = a^2 = 1$ yields $(x_1 x_2)^2 = x_1 x_2 x_1 x_2 = x_1^2 x_1 x_2^2 x_2 x_1 x_2 = [x_1, x_2] = a$. Thus, $x_1 x_2$ has order $4$ and so $D_3$ is non-singular. If $D_1$ and $D_2$ are singular, then $(x_1x_2)^2 = x_1 x_2 x_1 x_2 = [x_1, x_2] = a$. Thus, $D_3$ is non-singular.\qed

\medskip

Lemma \ref{one-or-all} provides us with a quadratic form $q$ with associated bilinear form $\beta$ on the dent space $\dents$, that is, $q(D_1 + D_2) = q(D_1) + q(D_2) + \beta(D_1, D_2)$. Namely, $q: \dents\map GF(2)$ which is given by

\begin{equation} \label{quadratic} q(D) = \left\{ \begin{array}{ll} 0 & \mbox{if } D \textnormal{ is singular} \\ 1 & \mbox{if } D \textnormal{ is non-singular} \end{array} \right.
\end{equation}

The forms $\beta$ and $q$ are called the \emph{bilinear} and \emph{quadratic forms associated} with the biextraspecial group $G$ and its dent space $\dents$.

\begin{defn} We define the \emph{type} of a biextraspecial group to be the type of its associated quadratic form $q$.
\end{defn}

Hence, a biextraspecial group has an even rank $m$ and type $\eps \in \{-,+\}$.

%------------ End of Section ---------------%

\section{Composition and decomposition}
\label{sec:decom}

Let $G = Q:L$ be a biextraspecial group with dent space $\dents$, and quadratic and bilinear forms $q$ and $\beta$ respectively. We adopt the language of vector spaces with respect to $\beta$ and $q$ in this and subsequent sections. Suppose that ${\cal U}$ is a nondegenerate subspace of $\dents$. Let $Q_1$ be the group generated by the dents of ${\cal U}$. Define the group \emph{associated} to ${\cal U}$ to be $G_1 = Q_1: L$. This semi-direct product is well defined because $L$ normalises each dent and so also $Q_1$.

Suppose that $\dents$ is the orthogonal sum of two nondegenerate subspaces ${\cal U}_1$ and ${\cal U}_2$ of dimension $k$ and $m-k$ with associated groups $G_1 = Q_1 : L$ and $G_2 = Q_2 : L$, respectively. The decomposition $\dents = {\cal U}_1 \perp {\cal U}_2$ implies that $Q$ can be written as the central product of $Q_1$ and $Q_2$. If $q\in Z(Q_1)$, then $[q, Q_1] = [q, Q_2] = 1$, and thus $[q, Q]=1$. Therefore $Z(Q_1) = Z(Q) = Z$. The quotient $Q_1/Z$ is a subspace of $Q/Z$ and thus has only natural modules as composition factors. Therefore $G_1$ is a biextraspecial group of rank $k$. It is clear that the dent space of $G_1$ is ${\cal U}_1$. Exactly the same argument can be used to show that $G_2$ is a biextraspecial group of rank $m-k$ with dent space ${\cal U}_2$. This argument applies to any finite number of orthogonal summands. So we have the following.

\begin{lem} \label{cordecom} Suppose that $\dents = {\cal U}_1 \perp\ldots\perp {\cal U}_k$ where each ${\cal U}_i$ is nondegenerate $n_i$-dimensional subspace. Let $G_i$ be the group associated to ${\cal U}_i$ for $i=1, \ldots, k$. Then $G_i$ is a biextraspecial group of rank $n_i$ with dent space ${\cal U}_i$.\qed
\end{lem}

The decomposition of $\dents$ into an orthogonal sum of nondegenerate subspaces translates to the decomposition of $G$ into smaller bi\-extraspecial groups. That is, given the decomposition in Lemma \ref{cordecom} we say that $G$ \emph{decomposes} into the biextraspecial groups $G_1$, \ldots, $G_k$. A question to ask is how can one assemble the smaller groups back together to produce the original biextraspecial group? We do this by introducing a composition operation between two biextraspecial groups.

Let $G_1 = Q_1 : L_1$ and $G_2 = Q_2 : L_2$ be biextraspecial groups of rank $m$ and $k$, respectively. Let $Z_i = Z(Q_i)$ for $i=1,2$. Let $\vphi: L_1\map L_2$ be an isomorphism and $\hat\vphi : Z_1 \map Z_2$ be the induced $L_1$-isomorphism. Let $\pi_i: G_i \map L_i$ be the projection map for $i=1,2$. Define $\hat G$ to be the subgroup of $G_1\times G_2$ consisting of the pairs $(g_1, g_2)$ such that $(g_1\pi_1)\vphi = g_2\pi_2$. Let $\hat Z$ be the subgroup of $Z_1\times Z_2$ consisting of the pairs $(z_1, z_2)$ such that $z_1\hat\vphi = z_2$. It is clear that $\hat Z$ is a normal subgroup of $\hat G$ since $Z_1$ and $Z_2$ are normal subgroups of $G_1$ and $G_2$, respectively. Let $L$ be the subgroup of $L_1\times L_2$ consisting of the pairs $(l_1, l_2)$ such that $l_2 = l_1\vphi$. Note that $Q_1\times Q_2$ and $L$ are subgroups of $\hat G$ and $\hat G = (Q_1\times Q_2):L$. Define $G = \hat G / \hat Z$ and $Q = (Q_1\times Q_2)/\hat Z$. Observe that $G$ is a split extension of $Q$ by $L$. Let $\hat Z(q_1, q_2)$ be in the centre of $Q$. Then $(q_1, q_2)$ commutes with $(q,1)$ for all $q\in Q_1$. That is, $(q_1,q_2)^{-1}(q,1)^{-1}(q_1, q_2)(q,1) = (q_1^{-1} q^{-1} q_1 q, 1) \in \hat Z$. This occurs if and only if $q_1$ commutes with all elements of $Q_1$, if and only if $q_1\in Z_1$. By symmetry, $q_2\in Z_2$. Therefore, $Z = Z(Q) = (Z_1\times Z_2)/\hat Z$. It is clear that $Z$ is the natural $L$-module since $Z_i$ is the natural $L_i$-module for $i=1,2$. Using the Third Isomorphism Theorem we obtain
\[ Q/Z = (Q_1\times Q_2)/\hat Z / (Z_1\times Z_2)/\hat Z \cong Q_1/Z_1 \times Q_2/ Z_2. \]
That is, $Q/Z$ is elementary abelian $2^{2(m+k)}$ and, in fact, $L$-isomorphic to $Q_1/Z_1 \times Q_2/Z_2$. Therefore, as $Q_1/Z_1$ and $Q_2/Z_2$ only have natural modules as composition factors, we have that $Q/Z$ only has natural $L$-modules as composition factors. In particular, $G$ is a biextraspecial group of rank $m+k$. The constructed group $G$ is denoted by $G_1\bistar G_2$ and is called the \emph{composition} of $G_1$ and $G_2$. We state this as a result.

\begin{thm} \label{thmcom} The group $G_1\bistar G_2$ is a biextraspecial group of rank $m+k$.\qed
\end{thm}

Note that this construction does not depend on the isomorphism between $L_1$ and $L_2$. Let $\vphi$ and $\psi$ be isomorphisms between $L_1$ and $L_2$. Then $\vphi\psi^{-1}\in\Aut(L_1) = \Inn(L_1)$, that is, there exists $g\in L_1$ such that $\psi = c_g \vphi$. The isomorphism $\psi$ induces the unique $L_1$-isomorphism  $\hat{c_g \vphi} = c_g \hat\vphi$. Let $G_\vphi = Q_\vphi : L_\vphi$ and $G_\psi = Q_\psi : L\psi$ be the composition of $G_1$ and $G_2$ with respect to $\vphi$ and $\psi$ respectively. Recall that $Q_\vphi = (Q_1\times Q_2)/\hat Z_\vphi$ and $Q_\psi = (Q_1\times Q_2)/\hat Z_\psi$, where $\hat Z_\vphi = \{ (z,z\hat\vphi) \mid z\in Z_1 \}$ and $\hat Z_\psi = \{ (z,(z^g)\hat\psi) \mid z\in Z_1\}$. The conjugation map $G_\vphi\map G_\psi$ given by $\hat Z_\vphi (q_1, q_2)(l, l\vphi) \mapto \hat Z_\psi(q_1, q_2^{g\vphi})(l,(l^g)\vphi)$ is an isomorphism.

\begin{rem} We can compose multiple biextraspecial groups without being concerned about the order of composition. That is, the composition $\bistar$ is associative and commutative up to isomorphism. In particular, given biextraspecial groups $G_1$, $G_2$, \ldots, $G_k$, we have that $G_1 \bistar G_2 \bistar \ldots \bistar G_k$ is uniquely defined up to isomorphism.
\end{rem}

We now generalise the composition operation to any finite number of biextraspecial groups.

\begin{cor} Let $G_1$, \ldots, $G_k$ be biextraspecial groups of rank $m_i$, \ldots, $m_k$, respectively. Then $G_1\bistar \ldots \bistar G_k$ is a biextraspecial group of rank $m_1+\ldots +m_k$.
\end{cor}

We show that the composition $\bistar$ is the reverse operation of decomposition. Let $G=Q:L$ be a biextraspecial group with dent space $\dents$. Let ${\cal U}_1$ and ${\cal U}_2$ be nondegenerate subspaces of $\dents$ with associated groups $G_1 = Q_1 : L$ and $G_2 = Q_2 : L$ such that $\dents = {\cal U}_1 \perp {\cal U}_2$. We take the identity isomorphism $L \map L$ which induces the identity $L$-isomorphism $Z\map Z$. We define $G_1\bistar G_2 = Q_{12} : L_{12}$ where $Q_{12} = (Q_1\times Q_2) / \hat Z$, $\hat Z = \{(z,z) \mid z\in Z\}$ and $L_{12} = \{ (l,l) \mid l\in L \}$. Note that $G = (Q_1\circ Q_2) : L$ since the orthogonality between ${\cal U}_1$ and ${\cal U}_2$ implies that every dent in $Q_1$ commutes with every dent in $Q_2$. Define the map $\Phi : G_1\bistar G_2 \map G$ by $\hat Z(q_1, q_2)(l,l) \mapto q_1 q_2 l$. A straightforward computation using the definition of multiplication and that $[Q_1, Q_2] = 1$ shows that $\Phi$ is an isomorphism. That is, $G$ and $G_1 \bistar G_2$ are isomorphic. This can be generalised to any finite number of biextraspecial groups.

\begin{lem} Let $G$ be decomposed into biextraspecial groups $G_1$, \ldots, $G_k$. Then $G$ and $G_1\bistar \ldots \bistar G_k$ are isomorphic.
\end{lem}

Let $\Fdents$ be the dent space of $G = G_1\bistar G_2$, and $\dents$ and $\Edents$ be the dent spaces for $G_1$ and $G_2$, respectively. The dents $D\in\dents$ and $E\in\Edents$ induce dents $\hat D = (D\times Z_2)/\hat Z$ and $\hat E = (Z_1\times E)/\hat Z$, respectively, of $\Fdents$.
In particular, using a simple counting argument we deduce that every dent of $\Fdents$ has the form $\lambda\hat D + \mu\hat E$ for $D\in\dents$ and $E\in\Edents$. This yields a natural isomorphism between $\Fdents$ and $\dents\oplus\Edents$, namely, that given by $\Psi:\lambda\hat D + \mu\hat E\mapto (\lambda D, \mu E)$. Let $q$, $q_1$, and $q_2$ be the quadratic forms associated with $\Fdents$, $\dents$, and $\Edents$, respectively. We define $(q_1,q_2)(D, E) = q_1(D) + q_2(E)$.  We make the following three observations. Firstly, $\hat D$ and $D$ have the same type and $\hat E$ and $E$ have the same type. In particular, $q(\hat D) = q_1(D)$ and $q(\hat E) = q_2(E)$. Secondly, $\hat D$ and $\hat E$ commute for all $D\in\dents$ and $E\in\Edents$. Finally, $\hat D_1$ and $\hat D_2$ commute in $\Fdents$ if and only if $D_1$ and $D_2$ commute in $\dents$. Similarly, $\hat E_1$ and $\hat E_2$ commute in $\Fdents$ if and only if $E_1$ and $E_2$ commute in $\Edents$. From these three observations we deduce that $q(\hat D + \hat E) = q_1(D) + q_2(E)$. As a result,
\[ (q_1,q_2)(\Psi(\lambda D + \mu E)) = (\lambda q_1(D), \mu q_2(E)) = \lambda q_1(D) + \mu q_2(E) = q(\lambda D + \mu E). \]
That is, $\Psi$ is an isometry between the orthogonal spaces $(\Fdents, q)$ and $(\dents, q_1)\oplus(\Edents, q_2)$. This can be generalised any finite number of dent spaces as follows.

\begin{lem} \label{lemcom} Let $G = G_1\bistar\ldots\bistar G_k$. Let $G$ have dent space $\dents$ with quadratic form $q$ and $G_i$ have dent space $\dents_i$ with quadratic form $q_i$ for $i=1,\ldots, k$. Then the orthogonal spaces $(\dents, q)$ and $(\dents_1, q_1)\oplus\ldots\oplus (\dents_k, q_k)$ are isometric.\qed
\end{lem}

We now state a standard result about orthogonal spaces. Given isometric spaces $(V,q)$ and $(U_1,q_1) \oplus\ldots\oplus (U_k, q_k)$ we have that if $\eps$ is the type of $q$ and $\eps_i$ is the type of $q_i$ for $i=1,2$, then $\eps = \eps_1 \ldots\eps_k$. Applying this to biextraspecial group we obtain the following result.

\begin{thm} Let $G$ be the composition of biextraspecial groups $G_1$, \ldots, $G_k$. If $G$ has type $\eps$ and $G_i$ has type $\eps_i$ for $i=1, \ldots, k$, then $\eps = \eps_1\ldots\eps_k$.\qed
\end{thm}

We have not shown that $\dents$ can be decomposed as the orthogonal sum of nondegenerate subspaces. However, this follows from the theory of orthogonal spaces. The orthogonal space $\dents$ can be decomposed as the orthogonal sum
\[ \dents = {\cal L}_1 \perp\ldots\perp {\cal L}_k\perp{\cal U},\]
where each ${\cal L}_i$ is a hyperbolic line and ${\cal U}$ is a $0$- or $2$-dimensional nondegenerate subspace. If $G$ is of positive type, then ${\cal U}$ has dimension $0$. If $G$ is of negative type, then ${\cal U}$ has dimension $2$. Each ${\cal L}_i$ corresponds to a biextraspecial group of rank $2$ of positive type and ${\cal U}$ (if nonzero) corresponds to a biextraspecial group of rank $2$ of negative type. We study and describe the implication of this in Section \ref{sec:class}.

%------------ End of Section ---------------%

\section{Extraspecial groups}
\label{sec:extrasp}
In this section we show that the dent space is isomorphic to an extraspecial $2$-group factored by its centre. By invoking the classification of extraspecial $2$-groups (which can be found in \cite{dh}) and a little extra work we can extract the properties of the dent space and the results it yields with this realisation.

For each involution $t\in L$ define $R_t = C_Q(t) = \{ q\in Q \mid q^t = q \}$ and $Z_t = Z(R_t)$. Let $Z = \gen{a,b}$ and $c=ab$. Throughout this section $t = t_a$ is a fixed involution which generates the stabiliser of $a$. We assume $L = \gen{s,t}$ such that $s$ sends $a$ to $b$ and $b$ to $c$. In particular, $t_b = ts^{-1}$ and $t_c = ts$ generate the stabilisers of $b$ and $c$, respectively. We first make the observation that $R_t^s = R_{t_b}$ and $R_t^{s^{-1}} = R_{t_c}$.

Recall that $\bar Q$ can be written as the direct sum of $m$ irreducible $2$-dimensional $L$-modules. Each such summand contains a $1$-dimensional subspace that is centralised by $t$. In particular, $C_{\bar Q}(t)$ is the direct sum of $m$ $1$-dimensional subspaces of $\bar Q$. 
If $q\in\bar Q$ is centralised by $t$, then $\bar q \bar q^s \bar q^{s^2} = \bar 1$. In particular, $\gen{\bar q, \bar q^s}$ is $2$-dimensional irreducible $L$-module of $\bar Q$. Furthermore $D = \gen{q, q^s, Z}$ is a normal subgroup of $G$. Therefore $D$ is a dent. We may assume that $x=q$ and $y=q^s$ is the standard basis for $D$ and so $q^t = q$. Therefore $t$ centralises exactly two elements of $\bar q = \{q,aq,bq,cq\}$, namely, $q$ and $aq$. That is, each element of $C_{\bar Q}(t)$ corresponds to exactly two elements of $R_t$, and so $R_t$ has cardinality $2^{m+1}$. For the rest of the section we assume that each dent has the standard basis.

\begin{lem} \label{extraspecial} The following assertions hold:
 \begin{itemize}
 \item[{\rm{(i)}}] $Q = R_t \times R_{t'}$ for any other involution $t'\in L$ not equal to $t$, and
 \item[{\rm{(ii)}}] $R_t$ is an extraspecial $2$-group of order $2^{m+1}$ with centre $Z_t = \gen{a}$.
 \end{itemize}
 \end{lem}

 \pf Let $t'$ be any other involution not equal to $t$, then either $R_{t_b} = R_{t'}$ or $R_{t_c} = R_{t'}$. Without loss of generality, suppose the former case holds and let $x\in R_t$. If $x^s\in R_t$, then $x^s = x^{st} = x^{s^{-1}}$ and in particular, $x^s = x$. But $s$ acts fixed-point-freely on $Q$. Hence $R_t \cap R_t^s = 1$. We know that $|R_t\times R_{t'}| = |R_t|\times |R_{t'}| = 2^{2 + 2m} = |Q|$, that is, $Q = R_t \times R_{t'}$. This proves (i). For (ii) note that $Z = Z_t \times Z_{t_b}$. Since $Z = \gen{a,b}$, $\{1, a\}\le Z_t$ and $\{1,b\}\le Z_{t_b}$, we have that $Z_t = \gen{a}$ and $Z_{t_b} = \gen{b}$. Let $x\in R_t$, then $(Z_tx)^2 = Z_t$ if and only if $x^2$ is contained in $Z_t$ and this follows since $x^2 \in Z\cap R_t = Z_t$. Therefore $R_t/Z_t$ is elementary abelian and thus $[R_t, R_t]\le\Phi(R_t) \le Z_t$. The derived subgroup $[R_t, R_t]$ is certainly not trivial and so equality holds.\qed

 \medskip
 For the next few results we let $R$ denote $R_t$.

 \begin{cor} \label{dentsextra} Every dent $D = \gen{x,y,Z}$ intersects $R$ in a group of order four. In particular,
\begin{itemize}
\item[{\rm{(i)}}] if $D$ is $2^4$, then $D\cap R = \gen{a, x} \cong 2^2$;
\item[{\rm{(ii)}}] if $D$ is $4^2$, then $D\cap R = \gen{x} \cong 4$.
\end{itemize}
\end{cor}

\pf Each of $\bar 1$ and $\bar x$ correspond to the subsets $\{1, a\}$ and $\{x, ax\}$ of $R$, respectively. Therefore $D\cap R = \{ 1, a, x, ax\}$. If $D$ is $2^4$, then $D\cap R = \gen{a,x}\cong 2^4$. If $D$ is $4^2$, then $D\cap R = \gen{x} \cong 4$.\qed

\medskip

Define $\bar R := R/Z_t$. Then $\bar R$ is elementary abelian of order $2^{2m}$ and can be viewed as an $m$-dimensional vector space over the field $Z_t = \{ 1, a\} \cong GF(2)$. Define the map $\beta_t: \bar R \times \bar R \map Z_t$ by $\beta_t(\bar x, \bar y) = [x, y]$ and $q_t: \bar R\map Z_t$ by $q_t(\bar x) = x^2$. It can verified that $q_t$ is a quadratic form with associated symmetric bilinear form $\beta_t$ and that $\beta_t$ is nondegenerate. Define $D_t = D \cap R$ for each dent $D$. Then $D_t$ is a group of order four as shown in the previous corollary. Each $\bar D_t = D_t/Z_t$ can be viewed as a $1$-dimensional subspace of $\bar R$. Define a map $\Psi:\dents\map\bar R$ by $D = \gen{x,y,Z}\mapto Z_t x$ and $0_\dents\mapto Z_t$. It is clear that $\Psi$ is bijective since it is injective and $\dents$ and $\bar R$ have the same cardinality. Let $D' = \gen{x', y', Z}$ and observe that $D + D' = \gen{ x x', y y', Z }$ where $xx'$ is centralised by $t$. Therefore $(D + D')\Psi = Z_t xx' = (Z_t x) (Z_t x') = (D\Psi)(D'\Psi)$. That is, $\Psi$ preserves the composition of dents. In particular, $\dents$ is a vector space and isomorphic to $\bar R$. Moreover, $q_t(D\Psi) = q_t(Z_t x) = x^2$. Therefore $q_t(D\Psi) = 0$ if and only if $x^2 = 1$, if and only if $D$ is $2^4$, if and only if $q(D) = 0$. This shows that $\Psi$ is an isometry.

\begin{lem} \label{dentsisoRt}The map $\Phi: \dents \map \bar R$ given by $D = \gen{x, y, Z} \mapto Z_t x$ is a linear isometry. That is, $(\dents, q)$ and $(\bar R, q_t)$ are isometric orthogonal spaces.\qed
\end{lem}

\medskip

The classification of extraspecial $2$-groups states that every extraspecial $2$-group $P$ has order $2^{2n+1}$ for some positive integer $n$ and that there are exactly two types of such a group for each $n$. The rank of $P$ is $n$ and the type of $P$ corresponds to the type of the quadratic form associated to it. In particular, the rank and type of $R$ corresponds to the type and rank of the biextraspecial group $G$. It also states that every extraspecial $2$-group can be written as the central product of extraspecial $2$-groups of smaller rank. Suppose that $R = R_1\circ\cdots\circ R_k$ where $R_i$ is an extraspecial $2$-group. Then $\bar R_i = R_i/Z_t$ is a nondegenerate subspace of $\bar R$ and $\bar R = \bar R_1 \perp \ldots \perp \bar R_k$ is an orthogonal sum and each $\bar R_i$ corresponds to a nondegenerate subspace ${\cal U}_i$ of $\dents$ via $\Psi$. In particular, this yields that $\dents$ can be written as $\dents = {\cal U}_1\perp\ldots\perp {\cal U}_k$. The classification of extraspecial $2$-groups more specifically states that $R$ can be decomposed into the central product of extraspecial groups of rank $2$ and this corresponds to the decomposition of $G$ into the product of biextraspecial groups of rank $2$.

%------------ End of Section ----------------------%

\section{Classification of biextraspecial groups}
\label{sec:class}

The problem of classifying biextraspecial groups is reduced to classifying the rank two case as mentioned at the end of Section \ref{sec:decom}. Let $G^+(2)$ and $G^-(2)$ be biextraspecial groups of rank $2$ of positive and negative type, respectively. Then the dent space of $G^+(2)$ is generated by two singular dents and the dent space of $G^-(2)$ is generated by two non-singular dents. We show that these two groups exist and are unique up to isomorphism.

Let $Q$ be the Sylow $2$-subgroup of $SL_3(4)$ that consists of uni-upper-triangular matrices. We identify $Q$ with the group of triples
\[\{(a,b,c) \mid a,b,c\in GF(4)\} \]
with binary operation $(a,b,c)(a',b',c') = (a+a', b+b', c + ab' + c)$ to make computations easier. Fix the subgroup $D = \{(a,a,b) \mid a,b\in GF(4)\}$ of order $16$ of $Q$ for the rest of this section. Let $Z$ denote the centre of $Q$ and observe that $Z = \{(0,0,a)\mid a\in GF(4)\} = [Q,Q]$. Let $\eta$ be a fixed element of $GF(4)$ of order $3$. Define an inner automorphism $s$ of $Q$ that is induced by the diagonal matrix $\rm{diag}(1,\eta, \eta^2)$. That is, $(a,b,c)^s = (\eta a, \eta b, \eta^2 c)$. Then $D$ and $Z$ are $\gen{s}$-invariant subgroups of $Q$. Observe that every element of $\bar Q := Q/Z$ has order $2$ and thus $\bar Q$ is elementary abelian $2^4$. In particular, $\bar Q$ is $\gen{s}$-module. By Maschke's theorem, $\bar Q$ is a completely reducible $\gen{s}$-module. The subspace $\bar D$ is a submodule and $s$ acts fixed-point-freely on $\bar Q$ which implies that $\bar Q$ decomposes as the sum of two $2$-dimensional irreducible submodules. The space $\bar Q$ can be viewed as a $2$-dimensional over the field $K$ of four elements where $K^* = \gen{s}$. We denote this space by $\bar Q_K$. The $1$-dimensional subspaces of $\bar Q_K$ correspond the $\gen{s}$-invariant subgroups of $Q$ of order $16$. In particular, there are five such subgroups. Let $U$ be such a $1$-dimensional space and $E$ be the corresponding subgroup of $Q$. Then $U = \{\bar 0, \bar x^s, \bar x^s, \bar x^{s^{-1}}\}$ and so $E = \gen{x, x^s, Z}$. In particular, $Z > [E,E] = \gen{[x,x^s]}$ is an $\gen{s}$-invariant subgroup of order $1$ or $2$. But $s$ acts fixed-point-freely on $Q$ and so $[E,E]=1$ and $E$ is abelian. Therefore, the five corresponding $\gen{s}$-invariant subgroups of order $16$ of $Q$ are abelian. Note that each pair of such subgroups intersect at $Z$. An element $(a,b,c)$ of $Q$ has order $2$ if and only if $(a,b,c)^2 = (0,0,ab) = 0$, if and only if $a=0$ or $b=0$. Therefore the two subgroups $D_1 = \{(a,0,b)\mid a,b\in GF(4)\}$ and $D_2 = \{(0,a,b)\mid a,b\in GF(4)\}$ are elementary abelian $2^4$ and $\gen{s}$-invariant. We know that $D$ is an $\gen{s}$-invariant subgroup of isomorphism type $4^2$. The two other subgroups, denoted by $D_3$ and $D_4$, must also be of isomorphism type $4^2$ since all elements outside of the $D_1\cup D_2$ have order $4$.

We now define two involutions of $Q$. Let $t$ be the automorphism induced by the unique non-trivial field isomorphism of $GF(4)$, that is, $(a,b,c)^t = (a^2, b^2, c^2)$. Let $t'$ be the automorphism given by the composition of $t$, the contragredient automorphism and the inner automorphism induced by the ``wrong diagonal'' matrix
\[ \left(\begin{array}{ccc}
    0& 0 & 1 \\
    0 & \eta & 0 \\
    \eta^2 & 0 & 0
    \end{array} \right).
\]
That is, $(a,b,c)^{t'} = (\eta^2 b^2, \eta^2 a^2, \eta(ab+c)^2)$. One can easily see that both $t$ and $t'$ are involutions. Set $L^+ = \gen{s,t}$ and $L^- = \gen{s, t'}$. Then $L^+\cong L^-\cong L_2(2)$. Define two semi-direct products $G^+ = Q:L^+$ and $G^- = Q:L^-$. Note that $Z$ is a natural $L^\eps$-module for $\eps\in\{0,1\}$. We show that $G^+$ and $G^-$ are biextraspecial groups of rank $2$ of positive and of negative type, respectively.

Note that $\bar Q$ is a $4$-dimensional $L^\eps$-module for $\eps\in\{-,+\}$. Furthermore, $\bar D$ is a $2$-dimensional irreducible $L^\eps$-submodule. The quotient $\bar Q / \bar D \cong Q/D$ is a $2$-dimensional $L^\eps$-module and can be shown to be irreducible. Using Lemma \ref{ext-VbyV}, we have that $\bar Q$ is the direct sum of two natural $L^\eps$-modules. The $L^\eps$-submodules of $\bar Q$ corresponds to $L^\eps$-invariant subgroups of order $16$ of $Q$. In particular, there are three $L^\eps$-invariant normal subgroups of $Q$ of order $16$ and they lie in the collection $\{D_1, D_2, D_3, D_4, D\}$. We know that $t$ normalises $D_1$, $D_2$ and $D$. The involution $t'$ normalises $D$, $D_3$ and $D_4$ since $t'$ interchanges $D_1$ and $D_2$. Therefore, $D_1$ and $D_2$ are normal subgroups of $G^+$, and $D_3$ and $D_4$ are normal subgroups of $G^-$. The group $Q$ can be generated by any two $D_i$ and in particular, $Q = \gen{D_1, D_2}$ and $Q = \gen{D_3, D_4}$. Thus $\bar Q = \bar D_1\oplus \bar D_2$ has only natural $L^+$-modules as composition factors and $\bar Q = \bar D_3 \oplus \bar D_4$ has only natural $L^-$-modules as composition factors. We conclude with the following theorem.

\begin{thm} The group $G^+$ is a biextraspecial group of rank $2$ of positive type with dent space $\gen{D_1, D_2}$ and $G^-$ is a biextraspecial group of rank $2$ of negative type with dent space $\gen{D_3, D_4}$.\qed
\end{thm}

We now wish to prove the uniqueness of these two types. Any two biextraspecial groups $G$ and $G'$ of rank $2$ and type $\eps$ have isometric dent spaces. In particular, there is a one-to-one correspondence between the dents of $G$ and $G'$ preserving the type of the dent. We use this fact to construct an isomorphism between $G$ and $G'$.

Let $G = Q:L$ and $G' = Q':L'$ be biextraspecial groups of rank $2$ and type $\eps$. Suppose that $Q = \gen{D_1, D_2}$ and $Q' = \gen{D_1', D_2'}$ and $Z = Z(Q)$ and $Z' = Z(Q')$ such that $D_i$ and $D_i'$ are isomorphic for $i=1,2$.  Let $\vphi: L\map L'$ be an isomorphism. We can define an isomorphism between $Q':L$ and $Q':L'$ by $ql \mapto q(l\vphi)$. In particular, we may assume that $L = L'$. Let $D_i$ and $D_i'$ have the standard basis for $i=1,2$ with respect the standard action of $L$. Thus there exists an isomorphism $\vphi_i: D_i\map D_i'$ such that $(d_i^l)\vphi_i = (d_i\vphi_i)^l$ for all $d_i\in D_i$ and for all $l\in L$ and that $\vphi_i$ sends $Z$ to $Z'$. In particular, $\vphi_i: Z\map Z'$ is an $L$-isomorphism. Furthermore, the maps $\vphi_1$ and $\vphi_2$ agree on $Z$. Indeed, let $z\in Z$ and $1\ne g\in L$ such that $z^g = z$. Then $(z^g\vphi_i) = z\vphi_i = (z\vphi_i)^g$. In particular, $z\vphi_1$ and $z\vphi_2$ are stabilised by $g$. Thus $z\vphi_1 = z\vphi_2$. Let $D_i = \gen{x_i, y_i, Z}$ and $D_i' = \gen{x_i', y_i', Z}$ have the standard bases. Define a map $\Phi: G \map G'$ by $(d_1 d_2 l)\Phi = (d_1\vphi_1) (d_2\vphi_2) l$.

\begin{thm} The map $\Phi$ is a group isomorphism.
\end{thm}
\pf We first show that $\Phi$ is a well defined. If $d_1 d_2= e_1 e_2$, then $e_1^{-1} d_1 d_2 e_2^{-1} = 1$. Then $(d_1d_2 l )\Phi = (e_1 e_2 l)\Phi$ if and only if $(e_1^{-1} d_1 d_2 e_2^{-1})\Phi = 1$, if and only if $\Phi$ is well defined on the identity of $G$. Let $d_1 d_2 = 1$. Then $d_1^{-1} = d_2 \in Z$. Therefore $(d_1 d_2) \Phi = (d_1\vphi_1)(d_1^{-1}\vphi_1) = 1$. This shows that $\Phi$ is well defined. It is clear that $\Phi$ is bijective since $\vphi_1, \vphi_2$ and $\vphi$ are bijective. It suffices to show that $\Phi$ is a homomorphism. Let $d_1 d_2 l, e_1 e_2 k \in G$. Then
\[ ((d_1 d_2 l)(e_1 e_2 k))\Phi = [d_2, e_1^l]\vphi_1[d_2\vphi_2, e_1^l\vphi_1] (d_1 d_2 l)\Phi (e_1 e_2 k)\Phi
\]
Set $e_1 = x_1^{m_1} y_1^{n_1}$ and $d_2 = x_2^{m_2} y_2^{n_2}$. Then, using the commutator table and that $\vphi: Z\map Z'$ is an $L$-isomorphism, we deduce that  $[d_2\vphi_2, (e_1^l)\vphi_1] = [d_2, e_1^l]\vphi_1$. Therefore $\Phi$ is a homomorphism.\qed

\medskip
This finishes the classification of biextraspecial groups. We state the classification in the following theorem.

\begin{thm} Let $G^\eps(m)$ denote a biextraspecial group of rank $m$ and type $\eps$. The $G^\eps(m)$ exists for each even positive integer $m$ and every type $\eps\in\{-,+\}$ and is unique up to isomorphism.
\end{thm}

\begin{rem}\label{out} We can generalise this argument for biextraspecial groups of any rank. Let $G$ and $G'$ be biextraspecial groups of rank $m$ and type $\eps$ with the same copy of $L_2(2)$. Let $Q = \gen{D_1, D_2}\circ\ldots\circ\gen{D_{m-1}, D_m}$ and $Q' = \gen{D_1', D_2'}\circ\ldots\circ\gen{D_{m-1}', D_m'}$ such that $D_i$ and $D_i'$ are of the same type. Then there exists an isomorphism $\vphi_i : D_i\map D_i'$ such that $(d_i^l)\vphi_i = (d_i\vphi_i)^l$ and that $\vphi_i$ maps $Z(Q)$ to $Z(Q')$. Define a map $\Phi: G\map G'$ by $(d_1\ldots d_m l)\Phi (d_1\vphi_1) \ldots (d_m \vphi_m) l$. Then $\Phi$ is a group isomorphism between $G$ and $G'$ that interchanges $D_i$ and $D_i'$ for all $i=1,\ldots, m$.\qed
\end{rem}

%------------ End of Section ------------------%

\section{The outer automorphism group}
\label{sec:out}
It is shown in this section that $\Out(G)$ acts linearly on the dent space while preserving the quadratic form $q$. Before we show this we establish a realisation of $\Out(G)$ in $\Aut(G)$.

\begin{lem} \label{Out_is_CAut} The outer automorphism group $\Out(G)$ is isomorphic to the group $C_{\Aut(G)}(L)$.
\end{lem}
\pf We identify $G$ with $\Inn(G)$ since $Z(G) = 1$. Lemma \ref{uniquecompl} states that there is exactly one class of complements to $Q$ in $G=Q:L$, that is, $G$ acts transitively on the set of complements to $Q$ in $G$. The Frattini Argument yields that $\Aut(G) = G \Aut(G)_L = G N_{\Aut(G)}(L)$. Using the Second Isomorphism Theorem and that $L$ is self-normalising in $G$, we have that $\Out(G)\cong N_{\Aut(G)}(L)/L$. Note that $Z(L) = 1$ and $\Aut(L) = \Inn(L) \cong L$ allowing us to identify $\Aut(L)$ with $L$ when convenient. Set $N=N_{\Aut(G)}(L)$ and define the map $\phi:N\map \Aut(L)$ by restriction, that is, $\alpha\mapto\alpha|_L$. Clearly this map is a homomorphism and has kernel $C_N(L)$. An element $\alpha\in \Aut(L)$ is a conjugation map via some element $x$ of $L$ and extends to an automorphism of $G$ that normalises $L$ by $\hat\alpha:ql\mapto (ql)^x = q^x l$. Thus $\phi$ is surjective. In particular, $N/C_N(L) \cong L$. But as $L$ is normal in $N$, and $C_N(L)$ and $L$ intersect trivially, we have that $N = C_N(L)\times L$. Thus $\Out(G) \cong C_N(L)L/L \cong C_N(L) = C_{\Aut(G)}(L)$.\qed

\medskip

Therefore, each element of $\Out(G)$ can be realised uniquely as an automorphism of $G$ that centralises $L$. The group $\Aut(G)$ is the split extension of $G = \Inn(G)$ by $C_{\Aut(G)}(L)$. We set $C := C_{\Aut(G)}(L)$ for the remainder of this section and study the action of $C$ on $\dents$.

\begin{lem} \label{out_acts_dents} The group $C$ acts linearly on the dent space preserving the quadratic form $q$.
\end{lem}

\pf Let $\alpha\in C$. Then $Z\alpha = Z$ as $Z$ is characteristic in $G$. Define an action of $L$ on $D\alpha/Z$ by $(Z(d\alpha))^l = (Zd^l)\alpha$. This induces an $L$-isomorphism between $D/Z$ and $D\alpha/Z$ since $\alpha$ commutes with the action of $L$ by definition. In particular, $D\alpha$ is a dent. Let $D_1$ and $D_2$ be dents. It can be verified that $(D_1+D_2)\alpha = D_1\alpha + D_2\alpha$. Furthermore, $D\cong D\alpha$, that is, $D$ and $D\alpha$ have the same type and thus $\alpha$ preserves the quadratic form $q$.\qed

\medskip
Thus, we have $\Out(G)$ acting on the dent space with an undetermined kernel and factor isomorphic to a subgroup of the orthogonal group. We now determine the kernel and show that this factor group is the full orthogonal group. Recall that Lemma \ref{uniqueauto} states that each dent has a unique automorphism that commutes with the action of $L$.

\begin{rem} \label{onedent} Consider the decomposition $Q = \gen{D_1, D_2}\circ\ldots\circ \gen{D_m, D_{m-1}}$. Using Remark $\ref{out}$ with $G = G'$, $\phi_j$ to be the identity for $j\ne i$ and $\vphi_i$ to be the unique non-trivial automorphism of $D_i$ that centralises $L$, we obtain an automorphism $\Phi_i$ of $G$ that centralises $L$ such that $\Phi_i$ centralises $D_j$ for $j\ne i$ and $\Phi_i$ restricted to $D_i$ is $\vphi_i$.
\end{rem}

Let $K$ be the kernel of the action of $C$ on $\dents$. Then $K$ consists of the automorphisms of $G$ that centralise the dent space and $L$.

\begin{lem} \label{linear_func} If $\alpha\in K$, then for each dent $D$, $\alpha|_D$ is either trivial or the unique non-trivial automorphism. Define a map $\hat\alpha:\dents\map GF(2)$ by

\[ D\mapto \left\{ \begin{array}{ll} 0 & \mbox{if } \alpha \mbox{ centralises } D \\ 1 & \mbox{ otherwise.} \end{array} \right. \]

\noindent Then $\hat\alpha$ is linear.
\end{lem}

\pf Let $D_1=\gen{x_1, y_1, Z}$ and $D_2=\gen{x_2, y_2, Z}$ be two dents with the standard bases. Suppose that $D_1\hat\alpha=D_2\hat\alpha$. Then $\alpha$ either centralises both $D_1$ and $D_2$, or does not centralise either dent. In the former case, it is clear that $\alpha$ centralises $D_1+D_2$, and in the latter case $(x_1x_2)\alpha = ax_1 ax_2 = x_1 x_2$ and similarly $(y_1y_2)\alpha = y_1y_2$, that is, $\alpha$ centralises $D_1+D_2$. Suppose that $D_1\hat\alpha\neq D_2\hat\alpha$, then $\alpha$ centralises exactly one of the dents say $D_1$. Clearly $(x_1x_2)\alpha = x_1 a x_2 \neq x_1 x_2$, that is, $\alpha$ does not centralise $D_1+D_2$. In particular, we have shown that $(D_1+D_2)\hat\alpha = D_1\hat\alpha + D_2\hat\alpha$ as required.\qed

\medskip
This lemma shows us that each automorphism from $K$ defines a linear functional $\dents\map GF(2)$, that is, an element of the dual of the dent space $\dents^*$. Recall that the basis of $\dents^*$ is $\omega_i: \dents\map GF(2) = \{0,1\}$ given by $D_j\mapto\delta_{ij}$ for $i=1,\ldots, m$. We now show that two $K$ and $\dents^*$ are isomorphic as groups.

\begin{thm} Let $K$ be as before. Then $K$ is isomorphic to the dual of the dent space via $\alpha\mapto\hat\alpha$.
\end{thm}

\pf Let us first verify that this is indeed a homomorphism. Let $\alpha_1, \alpha_2\in K$. We must show that $D(\hat{\alpha_1+\alpha_2}) = D\hat\alpha_1 + D\hat\alpha_2$. If $D(\hat{\alpha_1+\alpha_2}) = 0$, then $\alpha_1|_D = \alpha_2|_D$, thus it is clear that $D\hat\alpha_1 + D\hat\alpha_2 = 0$. If $D(\hat{\alpha_1 + \alpha_2}) = 1$, then exactly one $\alpha_i$ does not centralise $D$. In particular, $1=D\hat{\alpha_1\alpha_2} = D\hat\alpha_1 + D\hat\alpha_2$. Suppose that $D\hat\alpha_1 = D\hat\alpha_2$ for all dents $D$, then $d\alpha_1 = d\alpha_2$. In particular, $d\alpha_1 \alpha_2^{-1} = d$ for all $d\in D$ and $D\in\dents$, that is, $\alpha_1\alpha_2^{-1}$ centralises $Q$. As $\alpha_1\alpha_2^{-1}$ centralises $L$ (by definition), $\alpha_1\alpha_2^{-1}$ is the identity automorphism of $C$ and thus $\alpha_1 = \alpha_2$. Therefore $\alpha\mapto\hat\alpha$ is injective. Let $\omega_i$ be a basis element of $\dents^*$ and let $\Phi_i$ be defined as it is in Remark \ref{onedent}. In particular, $\Phi_i\in K$ and $\Phi_i\map\omega_i$. Therefore this mapping is bijective.\qed

\medskip

The group $C$ acts on the dent space $\dents$ with kernel $\dents^*$. The spaces $\dents^*$ and $\dents$ are isomorphic as modules. Therefore $C/\dents$ is a subgroup of the orthogonal group $O^{\eps}_m(2)$. Let $\dents = \gen{D_1, D_2}\perp\ldots\perp\gen{D_{m-1}, D_m}$. Let $\alpha$ be an orthogonal map of $\dents$, that is, $\alpha:\dents\map\dents$ is a linear isomorphism that sends singular dents to singular dents and non-singular dents to non-singular dents. In particular, $\dents = \gen{D_1\alpha, D_2\alpha}\perp\ldots\perp\gen{D_{m-1}\alpha, D_m\alpha}$. Let $G = (\gen{D_1, D_2}\circ\ldots\circ\gen{D_{m-1}, D_m}):L$ and $G = (\gen{D_1\alpha, D_2\alpha}\circ\ldots\circ\gen{D_{m-1}\alpha,D_m\alpha}):L$. Using Remark \ref{out}, there exists an automorphism of $G$ that centralises $L$ such that $D_i$ is mapped to $D_i\alpha$ for all $i=1, \ldots, m$. This automorphism corresponds to the orthogonal map $\alpha$ and therefore $C/\dents$ contains the entire orthogonal group.

\begin{thm} The group $C$ is an extension of the dent space $\dents$ by $GO^{\eps}_m(2)$.\qed
\end{thm}

Let $R = C_Q(t)$ as defined in Section \ref{sec:extrasp}. It is clear that $\Aut(R)$ acts on $\bar R = R/Z(R)$. Indeed, the inner automorphism group of $R$ is isomorphic to $R/Z_t = \bar R$ and a normal subgroup of $\Aut(R)$. Therefore $\Aut(R)$ acts on $\bar R$ by conjugation and the action is given by $(Z_t x)^\alpha = Z_t (x\alpha)$ for $\alpha\in\Aut(R)$. We now state a result that is analogous to Lemma \ref{out_acts_dents}.

\begin{lem} The group $\Aut(R)$ acts on $\bar R$ linearly preserving its quadratic form $q_t$.
\end{lem}

\pf It is quite clear that $\Aut(R)$ acts linearly on $\bar R$. Observe that $Z_t$ is a characteristic subgroup of order $2$ of $R_t$ and thus is centralised by $\Aut(R)$. In particular, $q_t(Z_t x\alpha) = (x\alpha)^2 = x^2\alpha = x^2 = q_t(Z_t x)$.\qed

\medskip
The last lemma is true for all extraspecial $2$-groups. We see in the next result that $\Aut(R)$ and $C$ are isomorphic.

\begin{lem} The group $C$ is isomorphic to the automorphism group of $R_t$ via a restriction map.
\end{lem}

\pf Let $\alpha\in C$ and let $x\in R_t$. Since $\alpha$ commutes with $L$, it commutes with $t$ and hence $x\alpha\in R_t$. Therefore $\alpha$ restricts to an automorphism of $R_t$. In particular, the restriction map $C\map \Aut(R_t)$ is a homomorphism. Let $s\in L$ be an element of order $3$ and recall that $Q = R_t\times R_t^s$. If $\alpha\in C$ centralises $R_t$, then it centralises $R_t^s$ and hence $Q$. Since $\alpha$ was chosen to centralise $L$, $\alpha$ is the trivial automorphism of $G$. Therefore $C\map \Aut(R_t)$ is injective and, furthermore, an isomorphism since both groups have the same cardinality.\qed

\medskip
We now identify $C$ with $\Out(G)$. Using Theorem 1 of \cite{g} we determine the nature of this extension, that is, whether it is split or not.

\begin{thm} Let $G = G^\eps(m)$ be the biextraspecial group of rank $m$ and type $\eps$. Then $\Out(G)$ is an extension of $\dents$ by $GO^{\eps}_m(2)$. If $m=2,4$ then the extension is split, otherwise it is non-split.\qed
\end{thm}

\end{document}